\documentclass[11pt,twoside]{article}
\usepackage{informat}
\usepackage{epsfig}

\usepackage{amsmath,amsfonts,amsthm,marvosym}
\usepackage{eucal,mathrsfs,dsfont}
\usepackage{graphicx,alltt}
\usepackage{paralist}
\usepackage{color,multicol}
\usepackage{lscape}
\usepackage{longtable}
\usepackage{array,tabularx,tabulary,booktabs} % Дополнительная работа с таблицами
\usepackage{longtable}  % Длинные таблицы
\usepackage{multirow} % Слияние строк в таблице
\usepackage{wrapfig} % Обтекание рисунков и таблиц текстом

\newcommand{\R}{\mathbb R}

\newcommand{\prt}{\partial}

\newcommand{\be}{\begin{equation}}
\newcommand{\ee}{\end{equation}}

\begin{document}
  \title{Parameter estimation in diffusion models with low regularity coefficients}
   \author{Dmytro Ivanenko, Rostyslav Pogorielov \\
      Taras Shevchenko National University of Kyiv, Kyiv, Ukraine \\
      Email:ida@univ.kiev.ua and rpogorielov95@gmail.com \\
      }
  \titleodd{Parameter estimation in diffusion models with low \dots }
  \authoreven{Dmytro Ivanenko, Rostyslav Pogorielov}
  \received{July 28, 2020}

  \abstract{The article considers parameter estimation constructing such as quasi-maximum likelyhood estimation and one step estimation in statistical models generated by solution of stochastic differential equation. It has been developed a software for parameter estimating and has been presented correspondent testing and comparing.}

  \maketitle

\section{Introduction}
It has been considered a models generated by diffusion processes that appear in a lot of areas such that financial processes [1], oceanography [2], physics [3], medicine [4] etc. The problem of parameter estimation in correspondent equations is an important problem that has lots of solutions.

In this article a discretely observed process that is a solution of a stochastic differential equation is considered

\begin{equation}\label{eq1}
d X_t = A\left(\alpha, X_t\right)dt + \sigma\left(\beta, X_t\right)d W_t,\ X_0=x_0,
\end{equation}
where Wiener process $W$ and process $X$ are real numbers, $\sigma$ and $A$ are known functions, $\alpha$ and $\beta$ are unknown parameters.
Observation is held with a constant step $h$. 

It has been proven in the article [5] existence of transition probability density decomposition in terms of Hermite polynomials. Using this decomposition it was constructed quasi-likelyhood function that is continuously differentiable in order to estimate quasi maximum likelihood estimator (QMLE) and compare it with other existing estimators. 

The structure of article is following.

In the first substantial part it has been presented general information about model and gives all necessary formulas that are used in parameter estimating such as Hermite polynomials, Hermite functions and transition probability density decomposition. Besides, in this section it has been presented quasi-likelyhood estimator and conditions on model items.

Next part presents other estimators such as conditional least square. They can be used as start for one step estimation by Newton-Raffson method, that requires initial value of parameter. 

Final section presents computational aspects of parameter estimating. Regularization and optimization methods have been described. It has been presented numerical experiments with comparison of algorithms and estimators.It has been also presented the description of the software that was developed in the process of solving this problem.

\subsection*{Acknowledgements} This research was partially supported by the Alexander von Humboldt Foundation within the Research Group Linkage Programme {\it Singular diffusions: analytic and stochastic approaches} between the University of Potsdam and the Institute of Mathematics of the National Academy of Sciences of Ukraine.

	\section{General information and explanation}
	It has been obtained in article [5] transition probability density decomposition by Hermite polynomials using parametrix method. It gave possibility to construct QML functional. Besides, in the cited article formulas for the coefficients of the approximating scheme of the diffusion process given by the equation (1) have been derived. Under the order of accuracy following notations [5] the following will be implied.
	
	It is said that a scalar function $f_{t}(\theta,x,y), t \in (0,1], x,y, \in \R$, is of at least order r if there exist some positive constants $C$ and $c$ which do not depend on $(\theta,x,y) \in\Theta\times\R$ and $t \in (0,1]$ such that
	$$
	|f_{t}(\theta,x,y)| \leq Ct^{r}\phi_{ct}(x-y)
	$$
	
	It is also assumed that drift and diffusion function satisfy local Lipschitz and linear growth conditions.
	
	The accuracy order of decomposition is 1 and corresponds to Milstein's scheme. Thus, in this article, test samples have been generated using the Milstein scheme, and the QML functional has been constructed by decomposing the first-order density.

	First of all, the necessary notation has been introduced in order to perform asymptotic expansions.
	
	Denote: $T$ is interval, $B = \sigma^{2}$, $n = \left[\frac{T}{h}\right]$ is number of observations, Gauss kernel $\phi_t\left(x\right)$ as:
	$$
	\phi_t\left(x\right)=\left(2\pi t\right)^{-1/2} e^{-x^2/\left(2t\right)}.
	$$
	
	Let $H_m\left(x, t\right),\ x > 0,\ t > 0$ --- Hermite polynomials, which are determined by the formula:
	$$
	\left(-1\right)^m { \frac{\partial^m}{ \partial x^m }} \phi_t\left(x\right)=H_m\left(x, t\right) \phi_t\left(x\right),\ m \ge 0.
	$$
	In particular
	\begin{multline}
		H_0(x;t)=1,\ H_1(x;t)=\frac{x}{t},\ H_2(x;t)=\frac{x^2}{t^2}-\frac{1}{t}, $$
	\end{multline}

	For the solution of a stochastic differential equation, the Milstein scheme is used, that has the form:
	\begin{multline}
		X_{i+1} = X_{i} + A(X_{i}, \alpha)h + 
		\sigma(X_{i}, \beta)W_{i} + \\
		0.5B(X_{i}, \beta)B'(X_{i}, \beta)(W_{i}^{2}-h).
	\end{multline}

	%After generating a sample of solutions to the equation \eqref{eq1}, a method of parameter estimator has to be chosen. 
	It is considered in this article conditional least squares, one step, scoring and quasi-likelyhood estimator. 
	
	Denote: $\theta = (\alpha. \beta)$ is the vector of uknown parameters. For creating quasi-likelyhood estimation a contrast function is used:
	$$
	QL_{n} (\theta; x_{1}, ... ,x_{n}) = \prod_{k = 1}^{n}p_{t}(\theta; x_{k-1},x_{k})
	$$
	where 
	\begin{multline}
		p_{t}(x,y) =  f^x_t(\theta;x,y)\Big\{1 + \\ +\frac{1}{4}t^2B(\beta,x)\prt_xb(\beta,x)H^x_3(\theta;x,y;t)\Big\}
	\end{multline}

	$$
	f^z_t(\theta;x,y)=\frac{1}{\sqrt{2\pi tb(\beta,z)}}\exp\Big\{-\frac{(y-x-A(\alpha,z)t)^2}{2B(\beta,z)t}\Big\},	
	$$	
	With respect to [1] it is assumed that drift and diffusion satisfy the following conditions:
	\begin{itemize}  \item[(i)] For some $\gamma>0$, one has 
		$A \in C_{B}^{0,1+\gamma};$ and $B = \sigma^{2} \in C_{B}^{1,1+\gamma}$
		\item[(ii)] $B$ is uniformly elliptic function. That is, $B(\theta,x) \geq c \geq 0, \theta \in \Theta, x \in {\rm I\!R} $
	\end{itemize}

	A lot of approaches for MLE estimation have been proposed when exact maximum likelhood is impossible. It can be obtained by numerical solution of Fokker-Plank equation using finite difference, finite element or Chebyshev collocation methods [6]. Another approaches contain Monte-Carlo Markov Chain methods and Metropolis-Hastings algorithm [7]. It can be also estimated by characterictic function [8] or indirect estimation [9]. Methods, proposed here require more computer resources and error is very difficult to control.
	
	An algorithm here is proposed based on decomposition from [5]. Similar approach was proposed by Ait-Sahalia in [10]. In this article it has been used other start approximation that led to over assumptions on the model and the algorithm of coefficients obtaining is too complicated.

	The quasi-maximum estimation is obtained by founding maximum of function
	 $\ln QL_{n} (\theta; x_{1}, ... ,x_{n})$ or minimum of $-\ln QL_{n} (\theta; x_{1}, ... ,x_{n})$. 
	
	\subsection*{Example}
\begin{figure}[htb]
		\begin{center}
			\epsfig{file=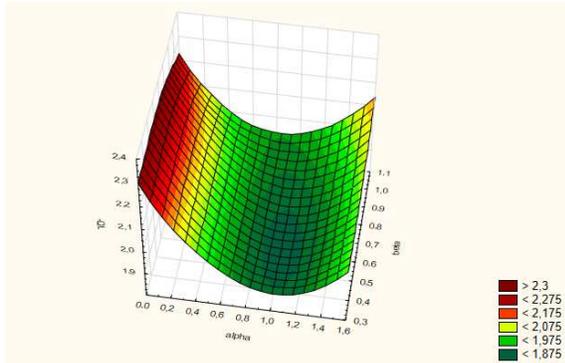,width=75mm} 
		\end{center}
		\caption{$-\ln QL_{n}$}
		\label{-QL_png}
	\end{figure}

	Figure \ref{-QL_png} shows minumum of log-likelyhood function in the case of $A = -x\alpha$, $B = 2 + \sin(x\times\beta)$ and $h = 0.5$, $T = 10000$, $\alpha = 1$, $\beta = 0.5$

	\section{Other estimators}
	
	\subsection{Conditional least squares}
	Conditional least square estimator is another estimation that can be found in given model. For discretization using Euler approximation the formula for loss function is given in [11] 
	$$Q(\alpha, \beta) =\sum_{k = 1}^{n}\dfrac{(X_{k}-X_{k-1} - A(x, \alpha)h)^{2}}{(B(x,\beta)^{2})h}$$
	
	In [12] there is proven weak and strong consistency of estimator as $T \rightarrow \infty$ and $T/n \rightarrow 0$. In [13] asymptotic normality and asymptotic effficiency are proven in case of known $\beta$
	
	In order to reduce bias it has been proposed conditional least squares estimator according to Milstein scheme. According to general definition of CLS estimator :
	$$Q(\alpha, \beta) =\sum_{k = 1}^{n}({X_{k}-g_{k}(\alpha,\beta, F_{k-1}))}^{2} $$
	where $g_{k}(\alpha,\beta, F_{k-1}) = E(X_{k}|F_{k-1})$ and $F_{k-1}$ is a sequence of of sub-sigma fields.
	
	For given discretely observed sample by Milstein approximation:
	
	\begin{multline}
	Q(\alpha, \beta) =\\
	\sum_{k = 1}^{n}\dfrac{(X_{k}-X_{k-1} - A(x, \alpha)h)^{2}}{E(B(X_{i})W_{i} + 0.5B(X_{i})B'(X_{i})(W_{i}^{2}-h))^{2}}
	\end{multline}

	Taking into consideration that:
	\begin{equation*}
	E[W^{p}] = 
	\begin{cases}
	0 &\text{p is odd}\\
	\sigma^{p} * (p-1)!! &\text{p is even}
	\end{cases}
	\end{equation*}
	and $W \sim N(0,\sqrt{h})$, it can be obtained:
	
	\begin{multline}
	Q(\alpha, \beta) = \\
	\sum_{k = 1}^{n}\dfrac{(X_{k}-X_{k-1} - A(x, \alpha)h)^{2}}{(B(X_{i}, \beta)^{2}h + 0.5B(X_{i}, \beta)^{2}B'(X_{i}^{2}, \beta)h^{2})}
	\end{multline}

	Advantage of CLS is in using less computational resourses but it is more biased than QMLE.
	
	Tests have shown that estimators LS and CLS - type strongly depend on the class of functions of the coefficients and in most cases simply do not work in the problems of estimating parameters when they are present in the diffusion coefficient. Thats why it has been proposed to use regularization functional of total type variation that has form $Q + |\beta - \beta_{0}|$
	
	\subsection*{Example}
\begin{figure}[!htb]
		\begin{center}
			\epsfig{file=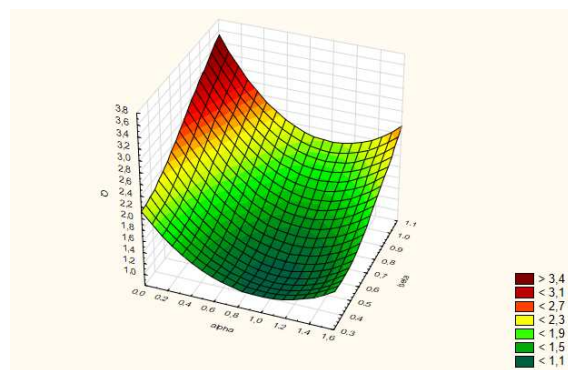,width=75mm} % figure1.eps is a standalone file in EPS format
		\end{center}
		\caption{Loss function with regularization}
		\label{QQ}
	\end{figure}

	Figure \ref{QQ} two shows loss function for functions $A = -x\alpha$, $B = 2 + \sin(x\times\beta)$ with parameters $h = 0.5$, $T = 10000$, $\alpha = 1$, $\beta = 0.5$, $\beta_{0} = 0.7$
	
	\subsection{One step and scoring estimators}
	As comparison of tables shows that Lp estimation is less accurate than QMLE. However, it can be used for one-step estimation, that asymptotically tend to MLE. The general formula for one-step estimation is:
	$$
	\theta_{one step} = \theta_{start} + H^{-1}(\theta_{start})*\nabla\ln QL_{n}(\theta_{start})
	$$ 
	where $\theta_{start}$ can be other estimator, for example CLS or QMLE, $H$ - Hesse matrix for $QL_{n}$ function, $\nabla$ - gradient vector.
	
	According to [5] partial derivatives have such form:
	\begin{multline}		
	\partial_{\alpha} p_{h}(\theta; x,y)=
	\Big\{1+\partial_{\alpha}A(\alpha,x)hH_{h}^{1,x}(\theta;x,y) \\ +\dfrac{h^{3}}{4}\partial_{\alpha}A(\alpha,x)\partial_{\beta}B(\beta,x)B(\beta,x)H_{h}^{4,x}(\theta;x,y) + \\
	(\dfrac{h^{2}}{4}\partial_{\alpha}A(\alpha,x)\partial_{\beta}B(\beta,x)) + \\ (\dfrac{h^{2}}{4}\partial_{\alpha x}A(\alpha,x)B(\beta, x))H_{h}^{2,x}(\theta;x,y)\Big\}f^x_t(\theta;x,y)
	\end{multline}
	
	\begin{multline}
	\partial_{\beta} p_{h}(\theta; x,y) = \Big\{\dfrac{h}{2}\partial_{\beta}B(\beta,x)H_{h}^{2,x}(\theta;x,y) + \\ \partial_{x}B(\beta,x)\partial_{\beta}B(\beta,x)B(\beta,x)H_{h}^{5,x}(\theta;x,y) + \\ \dfrac{h}{4}((B(\beta,x))\partial^{2}_{x\beta}B(\beta,x) + \\ \partial_{x}B(\beta,x)\partial_{\beta}B(\beta,x))H_{h}^{3,x}(\theta;x,y)\Big\}f^x_t(\theta;x,y)
	\end{multline}
	
	According to the rules of Hermite polynomials differentiating it can be obtained second derivatives for maximum likelyhood function:
	\begin{multline}
	\partial^{2}_{\alpha\alpha}p_t = 
	\Big\{h\partial^{2}_{\alpha\alpha}A(\alpha,x)H_{h}^{1,x} + t^{2}(\partial_{\alpha})^{2}A(\alpha,x)H_{h}^{2,x} \\ \dfrac{t^{3}}{4}\partial^{2}_{\alpha\alpha}A(\alpha,x)\partial_{x}B(\beta,x)B(\beta,x)H^{4,x} +
	\\
	\dfrac{t^{2}}{4}\partial^{2}_{\alpha\alpha}A(\alpha,x)\partial_{x}B(\beta,x)H^{2,x}+ 
	\\
	\dfrac{t^{2}}{2}\partial^{3}_{x\alpha\alpha}A(\alpha,x) B(\beta,x)H^{2,x}\Big\}f^x_t(\theta;x,y)
	\end{multline}
	\begin{multline}
	\partial^{2}_{\alpha\beta}p_t =
	\Big\{ 
	\dfrac{t^2}{4}\partial_{\alpha}A(\alpha,x)\partial_{\beta}B(\beta,x)H_{h}^{3,x}+\\
	\dfrac{t^3}{4}\partial_{\alpha}A(\alpha,x)\partial^{2}_{x\beta}B(\beta,x)B(\beta,x)H_{h}^{4,x}+
	\\
	\dfrac{t^3}{4}\partial_{\alpha}A(\alpha,x)\partial_{\beta}B(\beta,x)\partial_{x}B(\beta,x)H_{h}^{4,x}+\\
	\dfrac{t^4}{8}\partial_{\alpha}A(\alpha,x)\partial_{x}B(\beta,x)\partial_{\beta}B(\beta,x)B(\beta,x)H_{h}^{6,x}
	\\
	+ \dfrac{t^2}{4}\partial_{\alpha}A(\alpha,x)\partial^{2}_{x\beta}B(\beta,x)H_{h}^{2,x}+ 
	\\ \dfrac{t^3}{8}\partial_{\alpha}A(\alpha,x)\partial_{x}B(\beta,x)\partial_{\beta}B(\beta,x)H_{h}^{4,x}
	\\
	+
	\dfrac{t^2}{4}\partial^{2}_{x\alpha}A(\alpha,x)\partial_{beta}B(\beta,x)H_{h}^{2,x} +
	\\ \dfrac{t^3}{4}B(\beta,x)\partial^{2}_{x\alpha}A(\alpha,x)\partial_{beta}B(\beta,x)H_{h}^{4,x}\Big\}f^x_t(\theta;x,y)
	\end{multline}
	\begin{multline}
	\partial^{2}_{\beta\alpha}p_t =
	\Big\{ \dfrac{t^2}{4}\partial_{\beta}B(\beta,x)\partial_{\alpha}A(\alpha,x)H_{h}^{3,x} + \\ \dfrac{t^4}{8}B(\beta,x)\partial_{x}B(\beta,x)\partial_{\beta}B(\beta,x)H_{h}^{6,x} +
	\\
	\dfrac{t^3}{4}B(\beta,x)\partial^{2}_{x\beta}B(\beta,x)\partial_{\alpha}A(\alpha,x)H_{h}^{4,x} +\\ \dfrac{t^3}{4}\partial_{x}B(\beta,x)\partial_{\beta}B(\beta,x)\partial_{\alpha}A(\alpha,x)H_{h}^{4,x}\Big\}f^x_t(\theta;x,y)
	\end{multline}
	\begin{multline}
	\partial^{2}_{\beta\beta}p_t = 
	\Big\{\dfrac{t}{2}\partial^{2}_{\beta\beta}B(\beta,x)H_{h}^{2,x} + \dfrac{t^2}{4}(\partial_{\beta}B(\beta,x))^{2}H_{h}^{4,x} +
	\\
	\dfrac{t^3}{8}B(\beta,x)\partial_{x}B(\beta,x)\partial^{2}_{\beta\beta}B(\beta,x)H_{h}^{5,x} +\\ \dfrac{t^3}{8}B(\beta,x)\partial_{\beta}B(\beta,x)\partial^{2}_{x\beta}B(\beta,x)H_{h}^{5,x} +
	\\
	\dfrac{t^3}{8}(\partial_{\beta}B(\beta,x))^{2}\partial_{x}B(\beta,x)H_{h}^{5,x} + \\
	+ \dfrac{t^4}{16}B(\beta,x)(\partial_{\beta}B(\beta,x))^{2}\partial_{x}B(\beta,x)H_{h}^{7,x} 
	\\
	+ 
	\dfrac{t^2}{4}\partial_{\beta}B(\beta,x)\partial^{2}_{x\beta}B(\beta,x)H_{h}^{3,x} + \dfrac{t^2}{4}B(\beta,x)\partial^{3}_{x\beta\beta}B(\beta,x)H_{h}^{3,x} +
	\\
	\dfrac{t^3}{8}B(\beta,x)\partial^{2}_{x\beta}B(\beta,x)\partial_{\beta}B(\beta,x)H_{h}^{5,x} +\\ \dfrac{t^2}{4}\partial_{x}B(\beta,x)\partial^{2}_{\beta\beta}B(\beta,x)H_{h}^{3,x} + 
	\\
	+\dfrac{t^2}{4}\partial_{\beta}B(\beta,x)\partial^{2}_{x\beta}B(\beta,x)H_{h}^{3,x} + \\
	\dfrac{t^3}{8}(\partial_{\beta}B(\beta,x))^{2}\partial_{x}B(\beta,x)H_{h}^{5,x}\Big\}f^x_t(\theta;x,y)
	\end{multline}
	Denote ${h_{ij}}, i,j = 1,2$ the elements of Hessian. They are calculated via formulas:
	
	$h_{11} = (\dfrac{\partial_{\alpha}p_{t}}{p_{t}})^{2} - \partial^{2}_{\alpha\alpha}p_{t}$ 
	$h_{12} = \dfrac{\partial_{\alpha}p_{t}\partial_{\beta}p_{t}}{(p_{t})^{2}} - \partial^{2}_{\alpha\beta}p_{t}$
	\\
	$h_{21} = \dfrac{\partial_{\alpha}p_{t}\partial_{\beta}p_{t}}{(p_{t})^{2}} - \partial^{2}_{\beta\alpha}p_{t}$
	$h_{22} = (\dfrac{\partial_{\beta}p_{t}}{p_{t}})^{2} - \partial^{2}_{\beta\beta}p_{t}$
	
	Another way of obtaining estimation is scoring estimator where second derivatives are approximated by product of first derivatives.

	\section{Numerical experiments}
	
\begin{figure*}[!htb]
	
	\begin{center}
		\epsfig{file=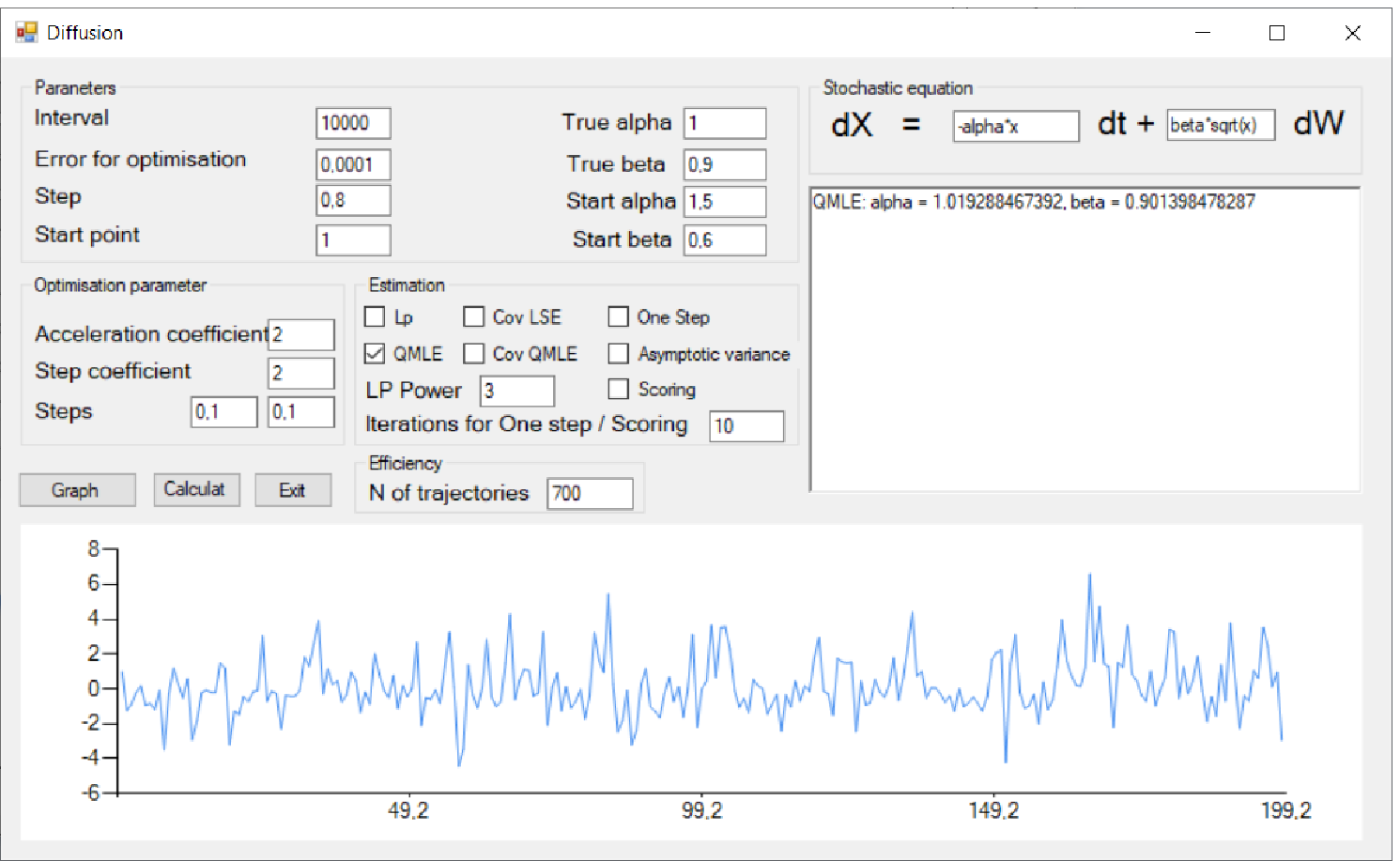,width=160mm} % figure1.eps is a standalone file in EPS format
	\end{center}
	\caption{Example for parameter estimation}
	\label{Sc}
\end{figure*}

	It has been developed a softmare for parameter estimation in diffusion model. It is possible to specify interval for process, step for discreet model, start point and true values of parameters $\alpha, \beta$ for Milstein simulation scheme. A wide variety of parameters is set for optimisation method.It has been chosen Hook-Jeeves method that uses pattern search. User can specify start point $(\alpha,\beta)$ as initial values for start iteration scheme. An example of parameter estimation is shown at \ref{Sc}

	Using software it has been obtained such results. Recall that it has been used Milstein scheme, regularization 
	with total variation type and Hook-Jeeves method. It has been chosen such values for optimization method:
	acceleration coefficient $= 1.1$, step coefficient $= 2$, steps $= 0.5$. 
	
	For $T = 10000$, $h = 0.8$, $A = -(\alpha*x)$, $B = 2+\sin(x*\beta)$, $\alpha = 1$, $\beta = 2$ ,$\alpha_{0} = 0.5$, $\beta_{0} = 1$ it has been calculated QMLE, CLS estimators and one-step and scoring estimators using QMLE and CLS estimators as start points. It has been simulated $1000$ trajectories and calculated mean and standart deviation of estimators. The results are in table 1.

	For $T = 10000$, $h = 0.8$, $A = \alpha(0.5-x)$, $B = \arctan(\beta*x)+2$, $\alpha = 1$, $\beta = 0.7$, $\alpha_{0} = 0.5$, $\beta_{0} = 0.5$ the results are in table 2.
	
	Tables show that QMLE estimator is more accurate than CLS and in most cases can be optained using optimization methods. Using one step or scoring estimators the result can be improved.
	
% bibliography style is not predefined; authors can apply their own,
% but it should be sent to matjaz.gams@ijs.si with the accepted paper

\begin{table*}[!htb]		
	
	\begin{center}
		
		\begin{tabulary}{\textwidth}{|c|c|c|c|c|}
			\hline 
			Estimation & Mean for $\alpha$ &  Standart deviation $\alpha$ & Mean for $\beta$ &  Standart deviation $\beta$  \\ 
			\hline 
			QMLE & 0,994644218 & 0,111534903 & 2,003946875 & 0,156352927\\ 
			\hline 
			CLS & 1,002137402 & 0,011515177 & 1,847436719 & 0,150998432 \\ 
			\hline 
			OS (QMLE) & 0,998466699 & 0,111517316 & 2,003347896 & 0,156133148 \\			
			\hline
			Scoring (QMLE) & 0,994661282 & 0,111525781 & 2,003346443 & 0,156187892 \\
			\hline
			OS (CLS) & 1,001578386 & 0,012638425 & 1,84981787 & 0,152179794 \\
			\hline
			Scoring (CLS) & 1,001657673 & 0,011502896 & 1,849437308 & 0,151560482 \\
			\hline
		\end{tabulary} 
		\caption{Comparison of estimators for $A = -(\alpha*x)$, $B = 2+\sin(x*\beta)$, $\alpha = 1$, $\beta = 2$ ,$\alpha_{0} = 0.5$, $\beta_{0} = 1$}
		\label{tab1}
		
	\end{center}
\end{table*}
\begin{table*}[!htb]		
	
	\begin{center}
		
		\begin{tabulary}{\textwidth}{|c|c|c|c|c|}
			\hline 
			Estimation & Mean for $\alpha$ &  Standart deviation $\alpha$ & Mean for $\beta$ &  Standart deviation $\beta$  \\ 
			\hline 
			QMLE & 1,02613125 & 0,009063776 & 0,737017188 & 0,016077046\\ 
			\hline 
			CLS & 1,000368359 & 0,007747116 & 0,578444141 & 0,017680205 \\ 
			\hline 
			OS (QMLE) & 1,001098021 & 0,00598422 & 0,720912942 & 0,073609449 \\			
			\hline
			Scoring (QMLE) & 1,009174357 & 0,004426152 & 0,691460678 & 0,019211411 \\
			\hline
			OS (CLS) & 1,000028765 & 0,10148162 & 0,650263072 & 0,133866962 \\
			\hline
			Scoring (CLS) & 1,000928822 & 0,027192407 & 0,741969292 & 0,111603981 \\
			\hline
		\end{tabulary} 
		\caption{Comparison of estimators for $A = \alpha(0.5-x)$, $B = \arctan(\beta*x)+2$, $\alpha = 1$, $\beta = 0.7$ , $\alpha_{0} = 0.5$, $\beta_{0} = 0.5$}
		\label{tab2}
		
	\end{center}
\end{table*}

\end{document}